\documentclass{amsart}
\usepackage{graphicx}
\usepackage{amscd}
\usepackage{amsmath}
\usepackage{amsfonts}
\usepackage{amssymb}
\theoremstyle{plain}

\newtheorem{corollary}{Corollary}

\newtheorem{definition}{Definition}
\newtheorem{example}{Example}

\newtheorem{proposition}{Proposition}
\newtheorem{remark}{Remark}

\newtheorem{theorem}{Theorem}
\numberwithin{equation}{section}

\begin{document}
\title[Reflexive Subspaces of $L_{1}$]{On P. Levy's Stable Laws and Reflexive Subspaces of $L_{1}$}
\author{Eugene Tokarev}
\address{B.E. Ukrecolan, 33-81, Iskrinskaya str., 61005, Kharkiv-5, Ukraine}
\email{tokarev@univer.kharkov.ua}
\dedicatory{Dedicated to the memory of S. Banach.}
\begin{abstract}To describe a set of functions, which forms a reflexive subspace $B$ of
$L_{1}$, a function $\eta(B,\varepsilon)$ that characterizes their average
integral growth is introduced: $\eta(B,\tau)=\underset{x\in B\backslash
\{0\}}{\sup}\{\underset{mes(e)=\tau}{\sup}\{\left\|  x(t)\chi_{e}(t)\right\|
_{L_{1}}/\left\|  x(t)\right\|  _{L_{1}}\}$.
\end{abstract}
\begin{abstract}It is shown that this function essentially depends on the geometry of $B$.
E.g., if $B$ contains a subspace isomorphic to $l_{p}$ for some $p<2$ then
$\eta(B,\tau)\leq c\tau^{1/q}$, where $q=p/(p-1)$. By the way, one question of
\ Ch. J. la Vall\'{e}e Poussin is answered and a short proof for a known
result on existence in every $L_{p}$ for $1<p<2$ of an uncomplemented subspace
isomorphic to the Hilbert space is obtained.
\end{abstract}
\begin{abstract}To describe a set of functions, which forms a reflexive subspace $B$ of
$L_{1}$, a function $\eta(B,\varepsilon)$ that characterizes their average
integral growth is introduced: $\eta(B,\tau)=\underset{x\in B\backslash
\{0\}}{\sup}\{\underset{mes(e)=\tau}{\sup}\{\left\|  x(t)\chi_{e}(t)\right\|
_{L_{1}}/\left\|  x(t)\right\|  _{L_{1}}\}$.
\end{abstract}
\begin{abstract}It is shown that this function essentially depends on the geometry of $B$.
E.g., if $B$ contains a subspace isomorphic to $l_{p}$ for some $p<2$ then
$\eta(B,\tau)\leq c\tau^{1/q}$, where $q=p/(p-1)$. By the way, one question of
\ Ch. J. la Vall\'{e}e Poussin is answered and a short proof for a known
result on existence in every $L_{p}$ for $1<p<2$ of an uncomplemented subspace
isomorphic to the Hilbert space is obtained.
\end{abstract}
\maketitle

\section{Introduction}

Consider the Lebesgue-Riesz space $L_{1}[0,1]=L_{1}$ of all summable real
functions defined on $[0,1]$. Certainly, $L_{1}$ is a Banach space under the
norm
\[
\left\|  x\right\|  _{L_{1}}=\int_{0}^{1}\left|  x\left(  t\right)  \right|
dt.
\]

Let $B$ be a (closed linear) subspace of $L_{1}$. Certainly, if $B$ contains a
disjoint sequence of functions, i.e., a such sequence $(f_{n})_{n<\infty
}=\left(  f_{n}\right)  $ that $f_{n}(x)f_{m}(x)=0$ (a.e.) when $n\neq m$,
then $B$ contains a subspace isometric to $l_{1}$ and, thus, is not reflexive.
The converse result was proved in [1].

\textit{Every non-reflexive subspace of }$L_{1}$\textit{\ contains an almost
disjoint in }$L_{1}$\textit{\ system of functions - i.e. such a sequence
}$\left(  g_{n}\right)  $\textit{\ that }$\underset{n\rightarrow\infty}{\lim
}\left\|  g_{n}-f_{n}\right\|  _{L_{1}}/\left\|  g_{n}\right\|  _{L_{1}}=0$
\textit{\ for some disjoint sequence }$\left(  f_{n}\right)  \subset L_{1}$.

This implies that any reflexive subspace $B$ of $L_{1}$ does not contain any
almost disjoint sequence.

In order to describe properties of a set of functions, which forms a reflexive
subspace $B$ of $L_{1}$, it will be convenient to introduce a
\textit{characteristic }$\eta$ \textit{of their average integral growth}.
Namely, let
\[
\eta(B,\varepsilon)=\sup\{\sup\{\left\|  x(t)\chi_{e}(t)\right\|  _{L_{1}%
}:\operatorname*{mes}(e)=\varepsilon\}:x\in B;\text{ \ }\left\|  x\right\|
_{L_{1}}=1\},
\]
where $\varepsilon>0$; $\operatorname*{mes}(e)$ is the Lebesgue measure of a
measurable subset $e\subset\lbrack0,1]$; $\chi_{e}\left(  t\right)  $ is
an\textit{\ indicator function }of $e$: $t\in e\Rightarrow\chi_{e}\left(
t\right)  =1$;$\ $ $t\notin e\Rightarrow\chi_{e}\left(  t\right)  =0$.

Clearly, $\eta(B,\varepsilon)$ is a non decreasing function; $\eta(B,0)=0$;
$\eta(B,1)=1$.

The aim of the article is to study the following question:

\textit{How the function }$\eta(B,\varepsilon)$\textit{\ depends on
geometrical properties of a given subspace }$B$\textit{\ of }$L_{1}$?

To state the main results, let us introduce for any non-decreasing function
$\psi\left(  t\right)  $ its \textit{involution} $\psi^{\star}\left(
t\right)  $, which is given by equalities:
\[
\psi^{\star}\left(  t\right)  =\sup\{h/\psi\left(  t\right)  :0<h\leq
t\}\text{ for \ }t\in(0,1];
\]%
\[
\psi^{\star}\left(  0\right)  =0.
\]

It will be shown that for any reflexive subspace $B$ of $L_{1}$, which
contains a subspace isomorphic to $l_{p}$ ($p\in\left(  1,2\right)  $) there
exists a constant $c\left(  B\right)  $ such that
\[
\eta^{\star}(B,t)\geq c\left(  B\right)  t^{1/p};
\]

If $B$ is of infinite dimension and contains a subspace, which is isomorphic
to the Hilbert space, then there exists a constant $c_{1}\left(  B\right)  $
such that
\[
\eta^{\star}(B,t)\geq c_{1}\left(  B\right)  \left(  \log_{2}\left(
t^{-1}+1\right)  \right)  ^{-1/2}.
\]

D. Aldous [2] showed that any infinite dimensional subspace of $L_{1}$
contains a subspace, which is isomorphic to some $l_{p}$ where $p\in
\lbrack1,2]$. So, the last of the estimates for $\eta^{\star}(B,t)$ covers all
infinite dimensional subspaces of $L_{1}$. Certainly, if $B\hookrightarrow
L_{1}$ (the last symbol denotes that $B$ is a subspace of $L_{1}$) is of
finite dimension then the estimates above need not to be valid. E.g., for
$\ n<\infty$ and $R^{(n)}$ be an $n$-dimensional subspace of $L_{1}$ spanned
by the first $n$ Rademacher functions, $\eta(R^{(n)},t)\asymp t$ for $t$ be
close to $0$, and $\eta^{\star}(R^{(n)},t)$ is discontinues at the point $t=0$.

To state these results it will be needed to regard the space $L_{1}$ as a
member of a class of symmetric Banach function spaces. Recall the definition.

\begin{definition}
A Banach space $E$ of (classes of) Lebesgue measurable real-valued functions
acting on the interval $[0,1]$ is called to be \textit{symmetric }if for any
(measurable) functions $x=x\left(  t\right)  $ and $y=y\left(  t\right)  $
following conditions hold:

1. If $x\in E$ and $\left|  y(t)\right|  \leq\left|  x(t)\right|  $ (a.e) then
$y\in E$ and $\left\|  y\right\|  _{E}\leq\left\|  x\right\|  _{E}$.

2. If $x\in E$ and functions $\left|  y(t)\right|  $ and $\left|  x(t)\right|
$ are equimeasurable then $y\in E$ and $\left\|  y\right\|  _{E}=\left\|
x\right\|  _{E}$.
\end{definition}

Let $\mathcal{S}$ be a class of all symmetric Banach function spaces (they in
the future will be referred as to symmetric spaces). This class contains
Lebesgue-Riesz spaces $L_{p}$; Orlicz spaces $L_{M}$, Lorentz spaces
$\Lambda\left(  \varphi\right)  $, Marcinkiewicz spaces $\mathsf{M}(\varphi)$
and $\mathsf{M}_{0}(\varphi)$ (the last ones will be defined below) and so on.

By the way, one question of Ch.J. la Vall\'{e}e Poussin, which will be
formulated below, will be answered. Besides these results, it will be obtained
an alternative simple proof of the following known result:

\textit{Any space} $L_{p}$ $(1<p<2)$ \textit{contains an uncomplemented
subspace, which is isomorphic to the Hilbert space.}

Some results of the article were announced in [3].

\section{Definitions and notations}

Let $E$ be a symmetric space. Its\textit{\ fundamental function} $\varphi
_{E}\left(  t\right)  $ is given by
\[
\varphi_{E}\left(  \tau\right)  =\left\|  \chi_{\lbrack0,\tau]}\left(
t\right)  \right\|  _{E}.
\]
To any measurable function $x(t)$ there correspond:

\textit{the distribution function}
\[
n_{x}\left(  s\right)  =\operatorname*{mes}\left(  \{t\in\lbrack0,1]:x\left(
t\right)  >s\}\right)  ;
\]

\textit{the non-increasing rearrangement}
\[
x^{\ast}\left(  t\right)  =\inf\{s\in\lbrack0,\infty):n_{\left|  x\right|
}\left(  s\right)  <t\}.
\]
With each symmetric space $E$ will be associated two more spaces:

\textit{the absolute continuous part}
\[
E_{0}=\{x\left(  t\right)  \in E:\underset{\varepsilon\rightarrow0}{\lim
}\left\|  x^{\ast}\left(  t\right)  \chi_{\lbrack0,\varepsilon]}\left(
t\right)  \right\|  _{E}=0\};
\]

\textit{the dual space}
\[
E^{\prime}=\{f(t)\in L_{1}:\int_{0}^{1}f(t)x(t)dt<\infty;\text{ \ }x\in
E;\text{ \ }\left\|  x\right\|  _{E}=1\}.
\]

The space $E_{0}$ will be equipped with a restriction of the norm $\left\|
\cdot\right\|  _{E}$; a norm of the space $E^{\prime}$ is given by
\[
\left\|  f\right\|  _{E^{\prime}}=\sup\{\left\|  f\cdot x\right\|  _{L_{1}%
}:x\in E;\text{ \ }\left\|  x\right\|  _{E}=1\}.
\]

Obviously, $\varphi_{E_{0}}(\tau)=\varphi_{E}(\tau)$; $\varphi_{E^{\prime}%
}(\tau)=\left(  \varphi_{E}(\tau)\right)  ^{\star}$.

Let $\mathcal{S}$ be a class of all symmetric spaces

The class $\mathcal{S}$ may be partially ordered by a relation $E\subset^{c}%
F$, which means that $E$, as a set of functions, is contained in $F$, and that
for any function $x\in E$, $\left\|  x\right\|  _{F}\leq c(E,F)\left\|
x\right\|  _{E}$. To be more precise, let us define an operator $i\left(
E,F\right)  :E\rightarrow F$, which asserts to any $x\in E$ the same function
$x\in F$. This operator is called \textit{the inclusion operator}. Its norm
(which is the infimum of all possible constants $c(E,F)$) is called
\textit{the inclusion constant}. So, the notation $E\subset^{c}F$ means that
the inclusion constant $c(E,F)<\infty$. If $c(E,F)=1$, then it will be written
$E\subset^{1}F$.

When $E\subset^{c}F$ and $F\subset^{c}E$, we shall write $E\approx^{c}F$.

In what follows it will be assumed that for any $E\in\mathcal{S}$ its
fundamental function $\varphi_{E}(t)$ satisfies the norming condition
$\varphi_{E}(1)=1$.\ This implies (see [4]) that a partially ordered set
$\left\langle E,\subset^{1}\right\rangle $ has both the maximal and the
minimal element: for any $E\in\mathcal{S}$,
\[
L_{\infty}\subset^{1}E\subset^{1}L_{1}.
\]

A function $\psi\left(  t\right)  $ defined on [0,1] will be called the
$\mathcal{M}$-\textit{function} provided

\begin{itemize}
\item $\psi\left(  0\right)  =0$; \ $\psi\left(  1\right)  =1$;

\item  Both $\psi\left(  t\right)  $ and $1/\psi\left(  t\right)  $ are
non-decreasing for $t>0$.
\end{itemize}

Let $\varphi\left(  t\right)  $ be the $\mathcal{M}$-function; $\varphi^{\sim
}\left(  t\right)  $ be the least concave majorant of $\varphi\left(
t\right)  $.

\textit{The Lorentz space }$\Lambda\left(  \varphi^{\sim}\right)  $ is a
Banach space of all measurable functions $x(t)$ on $[0,1]$ with the norm
\[
\left\|  x\right\|  _{\Lambda\left(  \varphi^{\sim}\right)  }=\int_{0}%
^{1}x^{\ast}\left(  t\right)  d\varphi^{\sim}\left(  t\right)  <\infty.
\]

\textit{The Marcinkiewicz space} $\mathsf{M}(\varphi)$ is a Banach space of
all measurable functions on [0,1] with the norm
\[
\left\|  x\right\|  _{\mathsf{M}(\varphi)}=\underset{0<h\leq1}{\sup}\{\frac
{h}{\varphi\left(  h\right)  }\int_{0}^{h}x^{\ast}\left(  t\right)
dt\}<\infty.
\]

The absolute continuous part $\left(  \mathsf{M}(\varphi)\right)  _{0}$ is
also called the Marcinkiewicz space and will be denoted by $\mathsf{M}%
_{0}(\varphi)$.

Clearly, $\varphi_{\Lambda\left(  \varphi^{\sim}\right)  }(\tau)=\varphi
^{\sim}\left(  \tau\right)  $; $\varphi_{\mathsf{M}\left(  \varphi\right)
}(\tau)=\left(  \varphi\left(  \tau\right)  \right)  ^{\star}$. Note that for
any $\mathcal{M}$-function $\varphi\left(  \tau\right)  $, its involution
$\ \varphi^{\star}\left(  \tau\right)  :=\left(  \varphi\left(  \tau\right)
\right)  ^{\star}$ is equal to $\tau/\varphi\left(  \tau\right)  $ for
$\tau>0$ and also is an $\mathcal{M}$-function.

It is known (see [4]) that spaces $\Lambda\left(  \varphi\right)  $ and
$\mathsf{M}(\varphi^{\star})$ are extreme among all symmetric spaces\ $E$ with
the concave fundamental function $\varphi_{E}\left(  \tau\right)
=\varphi\left(  \tau\right)  $ ($=\varphi^{\sim}\left(  \tau\right)  $); this
means that
\[
\Lambda\left(  \varphi\right)  \subset^{1}E\subset^{1}\mathsf{M}%
(\varphi^{\star}).
\]

Recall that spaces $\mathsf{M}_{0}(\varphi)$, $\Lambda\left(  \varphi\right)
$ and $\mathsf{M}(\varphi^{\star})$ are in duality:
\[
\left(  \mathsf{M}_{0}(\varphi)\right)  ^{\ast}=\Lambda\left(  \varphi\right)
\text{; \ }\left(  \Lambda\left(  \varphi\right)  \right)  ^{\ast}%
=\mathsf{M}(\varphi)\text{.}%
\]

\section{The Ch.J. la Vall\'{e}e Poussin's problem}

Let $B$ be a reflexive subspace of $L_{1}$; $E$ be a symmetric space. It will
be said that $B$ can be\textit{\ transferred }into $E$ if every function
$x(t)\in B$ belongs to $E$ and there exists a constant $c\in\lbrack1,\infty)$
such that for all $x\in B$
\[
c^{-1}\left\|  x\right\|  _{L_{1}}\leq\left\|  x\right\|  _{E}\leq c\left\|
x\right\|  _{L_{1}}.
\]

La Vall\'{e}e Poussin proved (cf. [5]), certainly, in other terminology, that
any reflexive subspace of $L_{1}$ can be transferred into some Orlicz space
$L_{M}$ distinct from $L_{1}$, and interested in bounds of such transfer. We
state the full solution of this problem.

Let $B$ be a reflexive subspace of $L_{1}$. Let $T(B)\subset\mathcal{S}$
denotes a class of all symmetric spaces $E$ with the property: $B$%
\textit{\ can be transferred into }$E$.

If $B$ is nonreflexive then the set $T(B)$ also can be defined. However in
this case it is trivial: if $E\in T(B)$ then $E$ is equal to $L_{1}$ up to
equivalent norm (more precise: $E\subset^{c}L_{1}$ and $L_{1}\subset^{1}E$,
i.e.,\ $E\approx^{c}L_{1}$). In what follows $B$ is assumed to be reflexive.

\begin{theorem}
For any reflexive subspace $B$ of $L_{1}$ there exists a symmetric space
$\frak{m}(B)$, which has following properties:

1. $B$ can be transferred into $\frak{m}(B)$;

2. For all $x\in B$ the following equality holds: $\left\|  x\right\|
_{L_{1}}=\left\|  x\right\|  _{\frak{m}(B)}$;

3. If $E\in T(B)$ then $\frak{m}(B)\subset^{c}E$.
\end{theorem}

\begin{proof}
Let $E\in T(B)$. Then there exists\ a constant $c<\infty$ (which depends on
both $B$ and $E$) such, that for all $x\in B$, $\left\|  x\right\|  _{E}\leq
c\left\|  x\right\|  _{L_{1}}$. Define on $E$ a new norm
\[
\left\|  \left|  x\right|  \right\|  _{E}=\max\{\left\|  x\right\|  _{L_{1}%
},\text{\ }c^{-1}\left\|  x\right\|  _{E}\}.
\]
Immediately, norms $\left|  \left\|  \cdot\right\|  \right|  $ and $\left\|
\cdot\right\|  _{E}$ are equivalent; $\left\langle E,\left|  \left\|
\cdot\right\|  \right|  \right\rangle =E_{1}$ is also a symmetric space;
$\varphi_{E_{1}}(1)=1$ and, for all $x\in B$, $\left\|  \left|  x\right|
\right\|  _{E}=\left\|  x\right\|  _{E}$. Consider, according to [6], a space
\[
\frak{m}(B)=\Cap\{E:E\in T(B)\},
\]
which consists of all functions $x(t)$ that are common to all $E\in T(B)$. A
norm on $\frak{m}(B)$ is given by
\[
\left\|  x\right\|  _{\frak{m}(B)}=\sup\{\left\|  \left|  x\right|  \right\|
_{E}:E\in T(B)\}.
\]
Obviously, $\left\|  x\right\|  _{\frak{m}(B)}=\left\|  x\right\|  _{L_{1}}$
for all $x\in B$ and $\frak{m}(B)\subset^{c}E$ for all $E\in T(B)$.
\end{proof}

\begin{theorem}
The fundamental function of the space $\frak{m}(B)$ is given by
\[
\varphi_{\frak{m}(B)}\left(  \tau\right)  =\eta^{\star}(B,\tau)=\underset
{0<h\leq\tau}{\sup}\{h/\eta\left(  B,\tau\right)  \}.
\]
\end{theorem}

\begin{proof}
Let $x\in B$; $\left\|  x\right\|  _{L_{1}}=1$. Then
\[
\int_{0}^{h}x^{\ast}(t)dt\leq\sup\{\int_{0}^{h}y^{\ast}(t)dt:y\in B;\text{
\ }\left\|  y\right\|  _{L_{1}}=1\}=\eta(B,h).
\]

Hence,
\[
\underset{0<h\leq1}{\sup}\eta^{-1}(B,h)\int_{0}^{h}x^{\ast}(t)dt=\left\|
x\right\|  _{\mathsf{M}(\psi)}\leq1,
\]
where $\mathsf{M}(\psi)$ is a Marcinkiewicz space with $\psi(t)=\left(
(\eta(B,t))^{\star}\right)  ^{\star}$. This implies that, for all $x\in B$,
$\left\|  x\right\|  _{\mathsf{M}(\psi)}=\left\|  x\right\|  _{L_{1}}$, i.e.
that $\mathsf{M}(\psi)$ belongs to $T(B)$. 

Clearly, $\frak{m}(B)\subset^{1}\mathsf{M}(\psi)$ and, hence,
\[
\varphi_{\frak{m}(B)}\left(  \tau\right)  \geq\eta^{\star}(B,\tau
)=\sup\{h/\eta\left(  B,\tau\right)  :0<h\leq\tau\}.
\]

To state the converse, note that from the H\"{o}lder inequality\ it follows
that for any $x\in L_{1}$
\begin{align*}
\left\|  x^{\ast}\chi_{\lbrack0,h]}\right\|  _{L_{1}} &  =\int_{0}^{h}x^{\ast
}(t)dt=\int_{0}^{1}[x^{\ast}(t)\chi_{\lbrack0,h]}(t)]\chi_{\lbrack0,h]}(t)dt\\
&  \leq\left\|  x^{\ast}\chi_{\lbrack0,h]}\right\|  _{E}\left\|  \chi
_{\lbrack0,h]}\right\|  _{E^{\prime}}=\varphi_{E^{\prime}}(t)\left\|  x^{\ast
}\chi_{\lbrack0,h]}\right\|  _{E},
\end{align*}
where $E$ is any symmetric space such that $x\in E$; $E^{\prime}$ is its dual.
Hence, for any $x\in B$,
\[
\varphi_{\frak{m}(B)}\left(  \tau\right)  \left\|  x^{\ast}\chi_{\lbrack
0,h]}\right\|  _{L_{1}}\leq h\left\|  x^{\ast}\chi_{\lbrack0,h]}\right\|
_{\frak{m}(B)}.
\]
Then
\begin{align*}
\varphi_{\frak{m}(B)}\left(  \tau\right)  \sup\{\left\|  x^{\ast}\chi
_{\lbrack0,h]}\right\|  _{L_{1}} &  :x\in B;\left\|  x\right\|  _{L_{1}%
}=1\}=\varphi_{\frak{m}(B)}\left(  \tau\right)  \eta(B,h)\\
&  \leq h\sup\{\left\|  x\right\|  _{\frak{m}(B)}:x\in B;\text{ \ }\left\|
x\right\|  _{L_{1}}=1\}=h.
\end{align*}
So, $\varphi_{\frak{m}(B)}\left(  \tau\right)  \leq\eta^{\star}(B,\tau)$.
\end{proof}

\section{How $\eta(B,t)$ depends on the geometry of $B$?}

We consider separately two cases:

\begin{enumerate}
\item $B$ contains a subspace, which is isomorphic to a space $l_{p}$ for some
$p\in(1,2)$.

\item $B$ contains a subspace, which is isomorphic to the Hilbert space
$l_{2}$.
\end{enumerate}

\begin{theorem}
Let $W$ be a reflexive subspace of $L_{1}$ that contains a subspace isomorphic
to $l_{p}$ for a given $p\in(1,2)$. Then
\[
\eta^{\star}(W,t)\geq c_{W}t^{1/p},
\]
where the constant $c_{W}$ depends only on $W$.
\end{theorem}

\begin{proof}
Let $p\in(1,2)$; $W$ be a subspace of $L_{1}$ isomorphic to $l_{p}$. Assume
that a sequence $\left(  f_{n}\right)  $ $(=\left(  f_{n}\right)  _{n<\infty
})\subset W$ is equivalent to the natural basis of $l_{p}$. It may be assumed
that $\left\|  f_{n}\right\|  =1$ for all $n\in\mathbb{N}$.

Fix $M\in\mathbb{N}$. Any function $f_{m}$ for $m<M$ generates a probability
measure $P_{m}$ on $[0,1]$, which is given by
\[
P_{m}(e)=\int\nolimits_{e}dn_{f_{m}}(s)\text{ \ for \ }e\subset\lbrack0,1].
\]
The set of functions $f_{1},$ $f_{2},$ $...,$ $f_{M}$ generates on a space
$[0,1]^{M}$ (the Carthesian product of $M$ copies of $[0,1]$) a probability
measure $P_{M}$ by the rule
\[
P_{M}\left(  e_{1}\otimes e_{2}\otimes...\otimes e_{M}\right)  =\Pi_{m=1}%
^{M}P_{m}\left(  e_{m}\right)  ,
\]
where $e_{i}\subset\lbrack0,1]$ for all $i\leq M$.

Let $\frak{G}_{n}$ be a group of rearrangements of the $M$-tuple
$\{1,2,...,M\}$.

Let a symmetric probability measure $\mu_{M}$ on $[0,1]^{M}$ is given by
\[
\mu_{M}\left(  e_{1}\otimes e_{2}\otimes...\otimes e_{M}\right)
=(M!)^{-1}\sum_{\sigma\in\frak{G}_{n}}P_{M}(e_{\sigma1}\otimes e_{\sigma
2}\otimes...\otimes e_{\sigma M}).
\]

For $M,N\in\mathbb{N}$, $N<M$, let $\mathsf{E}_{N}\mu_{M}$ be a symmetric
martingale of dimension $N$, which is associated with $\mu_{M}$ on $[0,1]^{M}$.

Let $D=D(M)$ be an ultrafilter on $\mathbb{N}$. Put
\[
\nu_{N}=\lim\nolimits_{D(M)}\mathsf{E}_{N}\mu_{M}%
\]

Since the set of probabilities $\{P_{M}:M\in\mathbb{N}\}$ (and, hence, both
$\{\mu_{M}:M\in\mathbb{N\}}$ and $\{\mathsf{E}_{N}\mu_{M}:M\in\mathbb{N\}}$)
is relatively weakly compact with respect to convergence by the law, $\nu_{N}$
is a probability measure.

Clearly, the ultrafilter $D$ defines an invariant mean on the union
$\cup\{\frak{G}_{n}:n<\infty\}$ and the family $\{\nu_{N}\}_{N<\infty}$ is
projective. So, it may be applied the classical Kolmogoroff theorem to define
in a standard way a probability measure $\nu_{\infty}$ on $[0,1]^{\mathbb{N}}%
$, which is consistent with $\{\nu_{N}\}_{N<\infty}$. Coordinate projections
$m:[0,1]^{\mathbb{N}}\rightarrow\lbrack0,1]$ generate a sequence of
independent random equidistributed functions $\left(  g_{m}\right)
_{m\in\mathbb{N}}$ (althou in a non-unique way: if $\tau:[0,1]\rightarrow
\lbrack0,1]$ is a preserving measure mapping, then functions $\{g_{m}%
(\tau^{-1}t):m\in\mathbb{N\}}$ also have same properties). Let $E$ be a
symmetric space in which $W$ can be transferred; $E\in T(W)$. By the symmetry
of $E$, $\left(  f_{n}\right)  \subset E$ implies that $\left(  g_{n}\right)
\subset E$. Of course,
\[
\left\|  \sum\nolimits_{n=1}^{M}g_{n}\right\|  _{E}\asymp M^{1/p}%
\]
because of $\left(  g_{n}\right)  $ is equivalent to the natural basis of
$l_{p}$. Hence, their mean distribution function converge to the standard
symmetric Levy's $p$-stable law. From [7] it follows that
\[
\lim\operatorname*{mes}\{t:n^{-1/p}\sum_{m=1}^{n}g_{n}(t)>\lambda
\}=C(p)\lambda^{-p},
\]
where
\[
C(p)=\left(  \int_{0}^{1}u^{-p}\sin(u)du\right)  ^{-1}.
\]

By standard arguments, similar to used in [8], the function $t^{-1/p}$ belongs
to the double dual $\left(  \frak{m}(W)\right)  ^{\prime\prime}$. 

This yields the inclusion $\frak{m}(W)\subset^{c}\mathsf{M}(t^{1/q})$, where
$q=p/(p-1)$. From the theorem 2 it follows that
\[
\eta^{\star}\left(  W,t\right)  =\varphi_{\frak{m}(W)}(t)\geq c_{W}%
\varphi_{\mathsf{M}(t^{1/q})}(t)=c_{W}t^{1/p}.
\]

So, the desired estimate is obtained.
\end{proof}

\begin{theorem}
Let $B$ be a reflexive infinite dimensional subspace of $L_{1}$. Then
\[
\eta^{\star}\left(  B,t\right)  \geq C_{B}[\log(1+t^{-1})]^{-1/2}%
\]
where $C_{B}$ depends only on $B$.
\end{theorem}

\begin{proof}
If $B$ contains a subspace isomorphic to some $l_{p}$, where $p\in(1,2)$, then
the desired estimate follows from Theorem 3. If $B$ does not contain any such
subspace, then, according to [2], $B$ contains a subspace $W$, which is
isomorphic to $l_{2}$. Let $\frak{m}(W)$ be the corresponding boundary space.
Obviously, $\frak{m}(W)\subset^{c}\frak{m}(B)$.

Clearly, $W$ (regarded as a subspace of $\frak{m}(W)$) does not contain any
almost disjoint in $\frak{m}(W)$ sequence. So, according to [9], $\frak{m}(W)$
satisfies the inclusion:
\[
G_{0}\subset^{c}\frak{m}(W),
\]
where $G_{0}$ is the absolutely continuous part of an Orlicz space $L_{M}$
with the norming function $M(u)$ equivalent to $\exp\left(  u^{2}\right)  -1$.
Its fundamental function up to a norming factor is equal to $[\log
(1+t^{-1})]^{-1/2}$. Hence,
\[
\eta^{\star}\left(  B,t\right)  =\varphi_{\frak{m}(B)}(t)\geq C_{1}%
\varphi_{\frak{m}(W)}(t)\geq C_{2}\varphi_{G_{0}}(t)\geq C_{B}[\log
(1+t^{-1})]^{-1/2}%
\]
for some constants $C_{1}$, $C_{2}$, $C_{B}$.
\end{proof}

\section{Boundary and binary symmetric spaces}

The following question is of interest:

\textit{What symmetric space }$E$\textit{\ may be regarded as the boundary
space }$\frak{m}(W)$ \textit{\ for some reflexive infinite dimensional
subspace }$W$ \textit{of }$L_{1}$\textit{? }

\begin{definition}
A symmetric space $E$ will be called \textit{a boundary space}, if there
exists such reflexive $B\hookrightarrow L_{1}$ (symbol $X\hookrightarrow Y$
means that $X$ is a subspace of $Y$) that $E\approx^{c}\frak{m}(B)$
\end{definition}

The class of all boundary spaces will be denoted by \textbf{Li}.

From [8] it follows that the above mentioned space $G_{0}$ is boundary for the
subspace $R\hookrightarrow L_{1}$, which is spanned by Rademacher functions.
From [10], th. 6, it may be deduced that the space \textsf{M}$_{0}\left(
t^{1/q}\right)  $ is the boundary space for the subspace $S_{p}\hookrightarrow
L_{1}$ that is spanned by a sequence $\left(  g_{n}(t)\right)  _{n<\infty}$ of
independent equidistributed functions with the characteristic function
\[
\overset{\sim}{g}\left(  t\right)  =\int_{0}^{1}\exp\left(  itg_{n}(u)\right)
du=\exp\left(  -\left|  t\right|  ^{p}\right)  ,
\]
where $p\in(1,2)$; $q=p/(p-1)$.

Nevertheless not every symmetric space is a boundary one. E.g., if
$G_{0}\nsubseteq^{c}E$, then each subspace of $E$ contains an almost disjoint
in $E$ sequence of functions and, hence, the corresponding norms $\left\|
\cdot\right\|  _{E}$ and $\left\|  \cdot\right\|  _{L_{1}}$cannot be
equivalent on infinite dimensional subspaces of $E$.

Another examples are given by members of a class of \textit{binary symmetric
spaces,} which was introduced in [11]. Let $E$ be a symmetric space;
$B\hookrightarrow E$;
\[
\eta_{E}(B,t)=\sup\{\left\|  x^{\ast}(t)\chi_{\lbrack0,\tau]}(t)\right\|
_{E}:x\in B;\text{ \ }\left\|  x\right\|  _{E}=1\};
\]%
\[
\eta_{E}\left(  B\right)  =\lim\nolimits_{t\rightarrow0}\eta_{E}(B,t)
\]

\begin{definition}
A symmetric space $E$ is said to be binary if for any infinite dimensional
subspace $B\hookrightarrow E$ its characteristic $\eta_{E}(B)$\ is equal
either to $0$ or $1$.
\end{definition}

The class of all binary spaces will be denoted by \textbf{Bi}.

The first example of a binary space was the space $L_{1}$. In [12] it was
shown that any Lorentz space is binary. In [11] this result was improved: it
was shown that any $L_{p}$ for $p<2$ is binary. For $p\geq2$ spaces $L_{p}$
are not binary ones.

\begin{theorem}
Classes of binary and of boundary symmetric spaces are not intersected:
$\mathbf{Li}\cap\mathbf{Bi=\varnothing}$.
\end{theorem}

\begin{proof}
If $B$ is a such subspace of $L_{1}$, which can be transferred into a binary
space $E$, then $\eta_{E}(B)=0$, as it follows from definitions. However, if
$E$ is the boundary space $\frak{m}(B)$ then $\eta_{E}(B)=1$ (and, hence, if
$E\approx^{c}\frak{m}(B)$ then $\eta_{E}(B)\neq0$). Indeed, from the proof of
theorem 2 it follows that for all $x\in B$
\[
\left\|  x\right\|  _{L_{1}}=\left\|  x\right\|  _{\frak{m}(B)}=\left\|
x\right\|  _{\mathsf{M}(\psi)},
\]
where $\psi(t)=\eta^{\star\star}(B,t)$. Hence,
\begin{align*}
1 &  \geq\eta_{\mathsf{M}(\psi)}(B,\tau)=\sup\{\left\|  x^{\ast}%
(t)\chi_{\lbrack0,\tau]}(t)\right\|  _{\mathsf{M}(\psi)}:x\in B,\left\|
x\right\|  _{L_{1}}=1\}\\
&  \geq\sup\{\underset{0<h\leq\tau}{\sup}(\eta(B,h))^{-1}\int_{0}^{h}x^{\ast
}(t)dt:x\in B,\left\|  x\right\|  _{L_{1}}=1\}\\
&  =\sup\{(\eta(B,\tau))^{-1}\underset{0<h\leq\tau}{\sup}\left\|  x^{\ast
}(t)\chi_{\lbrack0,h]}(t)\right\|  _{L_{1}}:x\in B,\left\|  x\right\|
_{L_{1}}=1\}=1
\end{align*}

Thus $\eta_{\mathsf{M}(\psi)}(B)=1$. Moreover,
\begin{align*}
1 &  \geq\eta_{\frak{m}(B)}(B,\tau)=\underset{\tau\rightarrow0}{\lim}%
\sup\left\|  x^{\ast}(t)\chi_{\lbrack0,h]}(t)\right\|  _{\frak{m}(B)}/\left\|
x\right\|  _{\frak{m}(B)}:x\in B\backslash\{0\}\}\\
&  \geq\underset{\tau\rightarrow0}{\lim}\sup\frac{\left\|  x^{\ast}%
(t)\chi_{\lbrack0,\tau]}(t)\right\|  _{\mathsf{M}(\psi)}}{\left\|  x\right\|
_{\mathsf{M}(\psi)}}\cdot\frac{\left\|  x\right\|  _{\mathsf{M}(\psi)}%
}{\left\|  x\right\|  _{\frak{m}(B)}}:x\in B\backslash\{0\}\}=\eta
_{\mathsf{M}(\psi)}(B)=1.
\end{align*}

These inequalities prove the theorem.
\end{proof}

The next result is, apparently, the first step in studying of order properties
of an order interval
\[
\mathsf{Kl}[\varphi]:=[\Lambda(\varphi);\mathsf{M}\left(  \varphi^{\star
}\right)  ]_{\subset^{1}}=\{E\in\mathcal{E}:\Lambda(\varphi)\subset
^{1}E\subset^{1}\mathsf{M}\left(  \varphi^{\star}\right)  \}
\]
of those symmetric spaces $E$, which have a given fundamental function
$\varphi_{E}(t)=\varphi\left(  t\right)  $ (which here is assumed to be
concave: if $\varphi\left(  t\right)  $ is not concave then the norm of the
corresponding Lorentz space $\Lambda(\varphi)$ does not satisfies the triangle inequality).

\begin{theorem}
Let $E,F\in\mathsf{Kl}[\varphi]$; $E$ be binary; $F$ be boundary. Then the
interval $[\Lambda(\varphi);E]_{\subset^{1}}$ contains no boundary spaces; the
interval $[F;\mathsf{M}\left(  \varphi^{\star}\right)  ]_{\subset^{1}}$
contains no binary spaces.
\end{theorem}

\begin{proof}
Let $H\in\lbrack F;\mathsf{M}\left(  \varphi^{\star}\right)  ]_{\subset^{1}}$.
Assume that $F=\frak{m}\left(  B\right)  $ for some $B\hookrightarrow L_{1}$.
By the preceding theorem, $\eta_{F}(B)=1$ and $\eta_{\mathsf{M}\left(
\varphi^{\star}\right)  }(B)=1$ as well. The same proof shows that $\eta
_{H}(B)=1$. If $F\approx^{c}\frak{m}\left(  B\right)  $ then the same argument
shows that $\eta_{H}(B)\neq0$. So, $H$ is not a binary space. If $H\in
\lbrack\Lambda(\varphi);E]_{\subset^{1}}$ and $H$ is a boundary space then
similarly $\eta_{H}(B)\neq0$ for some $B\hookrightarrow L_{1}$. Certainly this
implies that $\eta_{E}(B)\neq0$ what contradicts with $E\in\mathbf{Bi}$.
\end{proof}

\section{Uncomplemented Hilbert subspaces of $L_{p}$}

Let $Y$ be a Banach space. $X\hookrightarrow Y$ is said to be a
\textit{complemented subspace} of $Y$ if there exists a (bounded linear)
projection $P:X\rightarrow Y$.

M.I. Kadec and A. Pe\l czy\'{n}ski [1] showed that for $p\in\lbrack2,\infty)$
any subspace of $L_{p}$, which is isomorphic to the Hilbert space $l_{2}$, is
complemented in $L_{p}$.

For $p\in\lbrack1,2)$ the situation is different. Certainly, $L_{1}$ does not
contain any complemented reflexive subspace. For $p\in(1,2)$ the space $L_{p}$
contains a complemented subspace, which is isomorphic to the Hilbert space,
e.g., the space $R_{p}$ spanned in $L_{p}$ by Rademacher functions.

In [13], p. 145, it was pointed out that one W. Rudin's result on $\Lambda
(p)$-sets [14] implies that any space $L_{p}$ for $p\in\lbrack1,4/3]$ contains
an uncomplemented subspace $W$, which is isomorphic to $l_{2}$. It was
conjectured that the same is true for a whole interval $[1,2)$. This
conjecture was verified in [15], its proof given there used a quite
complicated construction. From the preceding technique a simple proof of this
fact follows.

\begin{theorem}
For any $p\in\lbrack1,2)$ the space $L_{p}$ contains an uncomplemented
subspace $W$, which is isomorphic to the Hilbert space.
\end{theorem}

\begin{proof}
Fix $p\in(1,2)$. Let $q=p/(p-1)$. Chose some $r\in(p,2)$. Let $s=r/(r-1)$. Let
$R$\ be a subspace of $L_{1}$, spanned by Rademacher functions. As it was
noted before, $\frak{m}(R)\approx^{c}G_{0}$, where $G$ is an Orlicz space
$L_{M}$ with norming function $M(u)=\exp\left(  u^{2}\right)  -1$;
$G_{0}=\left(  L_{M}\right)  _{0}$. It was also noted that
\[
\varphi_{G}\left(  \tau\right)  \asymp\varphi\left(  t\right)  :=[\log
(t^{-1}+1)]^{-1/2}%
\]
and, hence, $\eta\left(  R,t\right)  \asymp\psi\left(  t\right)
=t/\varphi\left(  t\right)  $ for $t>0$.

Chose a non-decreasing positive function $f(t)\in L_{1}$ such that
\[
\tau^{1/p}<\int_{0}^{\tau}f\left(  t\right)  d\psi\left(  t\right)
<\tau^{1/r}%
\]
and consider a subspace $W=[fR]$ of $L_{1}$, which consists of all measurable
functions $x(t)$ of kind: $x(t)=f(t)y(t)$ for some $y(t)\in R$. Clearly,
\[
\eta^{\star}(W,\tau)\asymp\tau(\int_{0}^{\tau}f(t)d\psi(t))^{-1}.
\]
Thus, since
\[
\varphi_{\frak{m}(W)}\asymp\eta^{\star}(W,\tau),
\]
there exists a constant $c=c(W)$ such that
\[
c^{-1}\tau^{1/s}<\varphi_{\frak{m}(W)}(\tau)<c\tau^{1/q},
\]

Hence, $W$ can be transferred into $L_{s}$ and cannot be transferred into
$L_{p}$. Let us regard the space $W$ as a subspace of $L_{s}$. This subspace
is isomorphic to $l_{2}$ and is complemented in $L_{s}$. Its conjugate,
$W^{\ast}$,\ which is also isomorphic to $l_{2}$, may be considered as a
complemented subspace of $\left(  L_{s}\right)  ^{\ast}=L_{r}$. Obviously, on
$W^{\ast}$ norms of $L_{r}$ and of $L_{1}$ are equivalent. So, there exists
some constant $\lambda<\infty$ such that for all $x\in W^{\ast}$ the following
inequality holds:
\[
\lambda^{-1}\left\|  x\right\|  _{L_{1}}\leq\left\|  x\right\|  _{L_{r}}%
\leq\left\|  x\right\|  _{L_{p}}\leq\left\|  x\right\|  _{L_{1}}\leq
\lambda\left\|  x\right\|  _{L_{p}},
\]
i.e., on $W^{\ast}$ norms of $L_{p}$ and of $L_{1}$ are equivalent as well.
$W^{\ast}$ may be regarded as a subspace of $L_{p}$, which is isomorphic to
$l_{2}$. Nevertheless $W^{\ast\ast}=W$ cannot be transferred into $L_{q}$.
This means that $W^{\ast}$ is an uncomplemented subspace of $L_{p}$. Since
$p\in(1,2)$ is arbitrary, the proof is done.
\end{proof}

\section{Appendix. On disjoint-like systems of functions}

Results of previous chapters show that disjoint sequences of functions play a
certain role in the theory of symmetric spaces: if a subset $K$ of a separable
symmetric space $E$ does not contain any almost disjoint system, then $K$ can
be transferred into $L_{1}$, i.e. on $K$ norms both of $E$ and of $L_{1}$ are
equivalent. However if $E$ is nonseparable, it may contain a such subset $K$
that consists of functions which is not almost disjoint althou which cannot be
transferred into $L_{1}$.

\begin{example}
Let a function $x(t)\in E\backslash E_{0}$. Let $\left(  e_{i}\right)  $ be a
system of measurable subsets of the $\operatorname*{supp}x\left(  t\right)
=\{t:x\left(  t\right)  \neq0\}$ with $e_{1}\supset e_{2}\supset...\supset
e_{i}\supset...$ , such that $\lim\operatorname{mes}e_{i}=0$.

Put $x_{n}\left(  t\right)  =x\left(  t\right)  \chi_{e_{n}}\left(  t\right)
.$ Then the set $K=\{x_{n}:n<\infty\}$ does not contain any almost disjoint
subsystem. However, on $K$ norms of $E$ and of $L_{1}$ are not equivalent:
since $\operatorname{mes}\operatorname*{supp}x_{n}$ tends to $0$, their
$L_{1}$-norm also has this property. And since $x(t)\in E\backslash E_{0}$,
there exists such $\varepsilon>0$ that $\left\|  x_{n}\right\|  _{E}%
\geq\varepsilon$ for all $n<\infty$.
\end{example}

Certainly, the system $\left(  x_{n}\left(  t\right)  \right)  $ may be called
disjoint-like. Here it will be presented some results on properties of
disjoint and disjoint-llike systems. For convenience recall some definitions.

\begin{definition}
Let $E$ be a symmetric space; $\{x_{n}\left(  t\right)  :n<\infty\}=\left(
x_{n}\right)  \subset E$. The system $\left(  x_{n}\right)  $ is said to be:

\begin{itemize}
\item  Disjoint, if $x_{n}\left(  t\right)  x_{m}\left(  t\right)  =0$ (a.e.)
when $m\neq n$.

\item  Almost disjoint in $E$ if there exists such disjoint system $\left(
y_{n}\right)  \subset E$, that
\[
\lim\nolimits_{n\rightarrow\infty}\left\|  x_{n}-y_{n}\right\|  _{E}/\left\|
x_{n}\right\|  _{E}=0.
\]

\item $\Delta_{E}$-system, if there exists such sequence of measurable subsets
$e_{n}\subset\left[  0,1\right]  $ ($n<\infty$) that
\[
\lim\nolimits_{n\rightarrow\infty}\operatorname{mes}e_{i}=1\text{ \ and
\ }\lim\nolimits_{n\rightarrow\infty}\left\|  x_{n}\chi_{e_{n}}\right\|
_{E}/\left\|  x_{n}\right\|  _{E}=0.
\]
\end{itemize}
\end{definition}

\begin{proposition}
If $E=E_{0}$ then every $\Delta_{E}$-system $\left(  x_{n}\right)  \subset E$
contains an almost disjoint in $E$ subsequence.
\end{proposition}

\begin{proof}
Since $E=E_{0}$, for every $x\in E$ and each $\varepsilon>0$ \ there exists
such real $v=v\left(  x,\varepsilon\right)  $ that $\operatorname{mes}e\leq v$
implies $\left\|  x\chi_{e}\right\|  _{E}\leq\varepsilon$.

Let $\varepsilon_{n}\searrow0$. Chose sequences $\left(  y_{n}\right)
\subset\left(  x_{n}\right)  $ and $\left(  e_{n}\right)  $; $e_{n}%
\subset\left[  0,1\right]  $ ($n<\infty$) such that $a_{n}=\operatorname{mes}%
e_{n}\nearrow1$; $\left\|  y_{n}\chi_{e_{e}}\right\|  \leq\varepsilon
_{n}\left\|  y_{n}\right\|  $ and
\[
\sum\nolimits_{n=j+1}^{\infty}\left(  1-a_{n}\right)  \leq v\left(
\varepsilon_{j},yj+1\right)  .
\]
Certainly this is possible.

Put $d_{n}=\left[  0,1\right]  \backslash e_{n}$; $D\left(  n\right)
=\cup_{m=n+1}^{\infty}d_{m}$. According to our choosing of $\left(
a_{n}\right)  $ and $\left(  y_{n}\right)  $
\[
\operatorname{mes}D\left(  n\right)  \leq1+a_{n}\text{; \ }\left\|  y_{n}%
\chi_{D\left(  n\right)  }\right\|  \varepsilon_{n}.
\]

Put $c_{j}=d_{j}\backslash D(j+1)$ ($j<\infty$). Sets $\left(  c_{j}\right)  $
are pairwice non-intersected and, hence, the system $\{\widehat{y}_{j}%
=y_{j}\chi_{c_{i}}\}$ is disjoint. So,
\begin{align*}
\left\|  y_{j}-\widehat{y}_{j}\right\|   &  =\left\|  y_{j}-y_{j}\chi_{d_{d}%
}+y_{j}\sum\nolimits_{j+1}^{\infty}\chi_{d_{d}}\right\|  \\
&  \leq\left\|  y_{j}-y_{j}\chi_{d_{d}}\right\|  +\left\|  y_{j}%
\sum\nolimits_{j+1}^{\infty}\chi_{d_{d}}\right\|  \\
&  \leq\left\|  y_{j}\chi_{e_{d}}\right\|  +\varepsilon_{j}\left\|
y_{j}\right\|  \leq2\varepsilon_{j}\left\|  y_{j}\right\|  .
\end{align*}
\end{proof}

\begin{remark}
If $E\neq E_{0}$ this theorem is not true as it follows from the example
before. So, it may be used as equivalent definition of the space $E_{0}$.
\end{remark}

To give a criterion to absence of $\Delta_{E}$-systems in a given subset $K$
of a symmetric space $E$ introduce some more notions.

Let $E\in\mathcal{S}$; $x\in E$; $K\subset E$. Put
\begin{align*}
S_{\varepsilon}^{E}\left(  x\right)   &  =\{t\in\left[  0,1\right]  :\left|
x\left(  t\right)  \right|  \geq\varepsilon\left\|  x\right\|  _{E}\};\\
M_{\varepsilon}\left(  E\right)   &  =\{x\in E:\operatorname{mes}%
S_{\varepsilon}^{E}\left(  x\right)  \geq\varepsilon\};
\end{align*}%
\[
\nu_{E}\left(  K\right)  =\lim\nolimits_{t\rightarrow0+}\inf\nolimits_{x\in
K\backslash\{0\}}\inf\nolimits_{\operatorname{mes}e\geq1-\varepsilon
}\{\left\|  x\chi_{e}\right\|  _{E}/\left\|  x\right\|  _{E}\}.
\]

\begin{theorem}
Let $E\in\mathcal{S}$; $K\subset E$. The following conditions are equivalent:

\begin{enumerate}
\item $K$ does not contain any $\Delta_{E}$-system;

\item  There exists a such $\varepsilon>0$ that $K\subset M_{\varepsilon
}\left(  E\right)  $;

\item $\nu_{E}\left(  K\right)  >0$.
\end{enumerate}
\end{theorem}

\begin{proof}
$\left(  1\Rightarrow3\right)  $. If $\nu_{E}\left(  K\right)  =0$ then there
are sequences $\left(  x_{n}\right)  \subset K$; $\left(  e_{n}\right)  $,
$e_{n}\subset\left[  0,1\right]  $ ($n<\infty$) and $\left(  \tau_{n}\right)
$, $\tau_{n}\searrow0$ such that
\[
\lim\nolimits_{n\rightarrow\infty}\left\|  x_{n}\chi_{e_{n}}\right\|
_{E}/\left\|  x_{n}\right\|  _{E}=0;\text{ \ }\operatorname{mes}e_{n}%
=1-\tau_{n}\nearrow1.
\]

Hence, $\left(  x_{n}\right)  $ is a $\Delta_{E}$-system.

$\left(  3\Rightarrow2\right)  $. Assume that\ $K\nsubseteq M_{\varepsilon
}\left(  E\right)  $ for any $\varepsilon>0$. Chose a sequence $\varepsilon
_{n}\searrow0$. For every $n<\infty$ there exists a function $x_{n}\in K$ such
that
\[
\operatorname{mes}\{t\in\left[  0,1\right]  :\left|  x\left(  t\right)
\right|  \geq\varepsilon\left\|  x\right\|  _{E}\}\leq\varepsilon_{n}.
\]

This yields that for a subset $K_{0}=\left(  x_{n}\right)  _{n<\infty}\subset
K$
\begin{align*}
\nu_{E}\left(  K_{0}\right)   &  \leq\lim\nolimits_{n\rightarrow\infty
}\left\|  x_{n}\chi_{\left[  0,1\right]  \backslash S_{\varepsilon}^{E}\left(
x_{n}\right)  }\right\|  _{E}/\left\|  x_{n}\right\|  _{E}\\
&  \leq\lim\nolimits_{n\rightarrow\infty}\varepsilon_{n}\left\|
x_{n}\right\|  _{E}/\left\|  x_{n}\right\|  _{E}=0.
\end{align*}

$\left(  2\Rightarrow1\right)  $. Let $K\subset M_{\varepsilon}\left(
E\right)  $ for some $\varepsilon>0$. Then for all $x\in K$%
\[
\operatorname{mes}\{t\in\left[  0,1\right]  :\left|  x\left(  t\right)
\right|  \geq\varepsilon\left\|  x\right\|  _{E}\}\geq\varepsilon.
\]

Hence, for every measurable subset $e\subset\left[  0,1\right]  $ of measure
$\operatorname{mes}e\geq1-\varepsilon/2$ and for all $x\in K$ there exists
such $d\subset e$ of measure $\operatorname{mes}d\geq\varepsilon/2$ that
$\left|  x\left(  t\right)  \right|  \geq\varepsilon\left\|  x\right\|  _{E}$
while $t\in d$. Thus,
\[
\left\|  x\chi_{e}\right\|  _{E}\geq\left\|  x\chi_{d}\right\|  _{E}%
\geq\left\|  \varepsilon\left\|  x\right\|  _{E}\chi_{d}\right\|
_{E}=\varepsilon\varphi_{E}\left(  \operatorname{mes}d\right)  \left\|
x\right\|  _{E}\geq\varepsilon\varphi_{E}\left(  \varepsilon/2\right)
\left\|  x\right\|  _{E},
\]
where $\varphi_{E}$ is a fundamental function of $E$.

Hence $K$ does not contain any $\Delta_{E}$-system.
\end{proof}

Subsets of symmetric spaces that does not contain any $\Delta_{E}$-system have
a series of important properties.

\begin{theorem}
Let $E\in\mathcal{S}$; $K\subset E$. The following conditions are equivalent:

\begin{enumerate}
\item $K$ does not contain any $\Delta_{E}$-system;

\item  The convergence by norm of $E$ is equivalent to the convergence by
measure on $K$.

If $E$ is differ from $L_{1}\left[  0,1\right]  $ (i.e. if $\lim
\nolimits_{t\rightarrow0+}\varphi_{E}\left(  t\right)  /t=0$) then both these
properties are equivalent to

\item  On $K$ norms $\left\|  \cdot\right\|  _{E}$ and $\left\|
\cdot\right\|  _{L_{l}}$ are equivalent.
\end{enumerate}
\end{theorem}

\begin{proof}
$\left(  1\Rightarrow2\right)  $. If $2$ is not true then there exists a
sequence $\left(  x_{n}\right)  \subset K$ such that $x_{n}/\left\|
x_{n}\right\|  $ convergents by measure to $0.$ Let $\varepsilon>0$. Put
\[
e_{n}=\{t\in\left[  0,1\right]  :\left|  x\left(  t\right)  \right|
\leq\varepsilon\left\|  x\right\|  _{E}\}.
\]

Then $\left\|  x\chi_{e}\right\|  _{E}/\left\|  x\right\|  _{E}\leq
\varepsilon$. Since $\operatorname{mes}e_{n}\rightarrow1$ (by our assumption),
$\nu_{E}\left(  K\right)  >0$ and, hence, $K$ contains some $\Delta_{E}$-system.

$\left(  2\Rightarrow1\right)  $. Assume that $1$ is not true. Then there
exists a sequence $\left(  x_{n}\right)  \subset K$ and measurable sets
$e_{n}\subset\left[  0,1\right]  $, $n<\infty$ with $\operatorname{mes}%
e_{n}\nearrow1$ such that
\[
\lim\nolimits_{n\rightarrow\infty}\left\|  x_{n}\chi_{e_{n}}\right\|
_{E}/\left\|  x_{n}\right\|  _{E}=0.
\]

This yields that $x_{n}\chi_{e_{n}}$ convergents by measure to $0$. Certainly,
the same is true for $\left(  x_{n}\right)  $.

$\left(  1\Rightarrow3\right)  $. According to the previous theorem, $1$
implies that $K\in M_{\varepsilon}\left(  E\right)  $ for some $\varepsilon
>0$. Hence for all $x\in K$%
\[
\left\|  x\right\|  _{E}\geq\left\|  x\right\|  _{L_{1}}\geq\left\|
x\chi_{S_{\varepsilon}^{E}\left(  x\right)  }\right\|  _{l_{1}}\geq\left\|
\varepsilon\left\|  x\right\|  _{E}\chi_{S_{\varepsilon}^{E}\left(  x\right)
}\right\|  _{L_{1}}\geq\varepsilon\varphi_{E}\left(  \varepsilon\right)
\left\|  x\right\|  _{E}.
\]

Now assume that $\lim\nolimits_{t\rightarrow0+}\varphi_{E}\left(  t\right)
/t=0$.

$\left(  3\Rightarrow1\right)  $. Let $1$ be not true. Then for every
$\varepsilon>0$ there exists such $x_{\varepsilon}\in K$ that
$\operatorname*{mes}S_{\varepsilon}^{E}\left(  x\right)  \leq\varepsilon$.

Chose some $\varepsilon_{0}>0$. Let $\varepsilon$ be such that $\varphi
_{E^{\prime}}\left(  \varepsilon\right)  +\varepsilon<\varepsilon_{0}$, where
$\varphi_{E^{\prime}}\left(  \varepsilon\right)  =\left(  \varphi_{E}\left(
\varepsilon\right)  \right)  ^{\ast}$ be the fundamental function of the dual
space $E^{\prime}.$

From the H\"{o}lder inequality it follows that%
\begin{align*}
\left\|  x_{\varepsilon}\right\|  _{L_{1}} &  =\int\nolimits_{S_{\varepsilon
}^{E}\left(  x_{\varepsilon}\right)  }\left|  x_{\varepsilon}\left(  t\right)
\right|  dt+\int\nolimits_{\left[  0,1\right]  \backslash S_{\varepsilon}%
^{E}\left(  x_{\varepsilon}\right)  }\left|  x_{\varepsilon}\left(  t\right)
\right|  dt\\
&  \leq\left\|  x_{\varepsilon}\right\|  _{E}\varphi_{E^{\prime}}\left(
\varepsilon\right)  +\varepsilon\left\|  x_{\varepsilon}\right\|  _{E}%
\leq\varepsilon_{0}\left\|  x_{\varepsilon}\right\|  _{E}.
\end{align*}

Since $\varepsilon_{0}$ is arbitrary, this inequality contradicts to the
equivalence of norms of $E$ and of $L_{1}$ on $K$.
\end{proof}

Before it was introduced the characteristic $\eta_{E}\left(  K\right)  $ of a
subset $K\subset E$. It also is related to question of extracting
disjoint-like sequences.

\begin{theorem}
If $\eta_{E}\left(  K\right)  <1$ then $K$ does not contain any $\Delta_{E}$-system.
\end{theorem}

\begin{proof}
It follows from the triangle inequality that$\left\|  x\chi_{S_{\varepsilon
}^{E}\left(  x\right)  }\right\|  _{E}\geq\left(  1-\varepsilon\right)
\left\|  x\right\|  _{E}$.

Therefore if $x\notin M_{\varepsilon}\left(  E\right)  $ then
\[
\sup\{\left\|  x\chi_{e}\right\|  _{E}:\operatorname*{mes}e\leq\varepsilon
\}\geq\left(  1-\varepsilon\right)  \left\|  x\right\|  _{E}.
\]

So, if $K\nsubseteq M_{\varepsilon}\left(  E\right)  $ then there exists such
$x_{\varepsilon}\in K$ that $\eta_{E}\left(  x_{\varepsilon},\varepsilon
\right)  \geq1-\varepsilon$.

Since $\varepsilon$ is arbitrary, $\eta_{E}\left(  K\right)  =1$.
\end{proof}

\begin{remark}
The converse is not true. Indeed, if \textit{symmetric separable space }%
$E$\textit{\ is a boundary space }$\frak{m}(W)$ \textit{for some reflexive
infinite dimensional subspace }$W$ \textit{of }$L_{1}$, then, as it was shown
before, $\eta_{E}\left(  W\right)  =1$ but $W$ does not contain any almost
disjoint sequence.

However, in some cases the converse result is true.
\end{remark}

\begin{theorem}
Let $E\in\mathcal{S}$; $E=E_{0}$ and $K\subset E$ be conditionally compact.
Then $\eta_{E}\left(  K\right)  =0$.
\end{theorem}

\begin{proof}
Let $\varepsilon>0$; $\{x_{1},x_{2},...,x_{n}\}$ be an $\varepsilon$-net for
$K$. Since $E=E_{0}$, there exists $\delta>0$ such that
\[
\sup\{\left\|  x_{i}\chi_{e}\right\|  _{E}:\operatorname*{mes}e\leq
\delta\}\leq\varepsilon/2\min\{\left\|  x_{i}\right\|  _{E}:i=1,2,...,n\}.
\]

From the triangle inequality it follows that \ $\left\|  x\chi_{e}\right\|
_{E}\leq$ $\varepsilon$ $\left\|  x\right\|  _{E}$ for all $x\in K$ whenever
$\operatorname*{mes}e\leq\delta$. Since $\varepsilon$ is arbitrary, $\eta
_{E}\left(  K\right)  =0$.
\end{proof}

\begin{corollary}
If $E=E_{0}$ and $B\hookrightarrow E$ is of finite dimension then $\eta
_{E}\left(  B\right)  =0$.
\end{corollary}

On every symmetric space $E$ it may be introduced an equivalent norm, say,
$\left|  \left\|  \cdot\right\|  \right|  $, such that for the renormed space
$\widetilde{E}=\left\langle E,\left|  \left\|  \cdot\right\|  \right|
\right\rangle $ the condition $\eta_{\widetilde{E}}\left(  K\right)  <1$ would
be equivalent to $\nu_{E}\left(  K\right)  >0$ (or any other of equivalent
conditions of theorems 8 and 9).

\begin{proposition}
Put $\left|  \left\|  x\right\|  \right|  =(\left\|  x\right\|  _{E}+\left\|
x\right\|  _{L_{1}})/2$. Then

\begin{enumerate}
\item  This norm is equivalent to $\left\|  \cdot\right\|  _{E}$;

\item  For every $K\subset E$ that does not contain any $\Delta_{E}$-system
$\eta_{\widetilde{E}}\left(  K\right)  <1$.
\end{enumerate}
\end{proposition}

\begin{proof}
1. $\left\|  x\right\|  _{E}\leq\left\|  x\right\|  _{E}+\left\|  x\right\|
_{L_{1}}\leq2\left\|  x\right\|  _{E}$.

2. According to the theorem 9 it may be assumed that on $K$ norms of $E$ and
of $L_{1}$ are equivalent. Clearly,
\[
\eta_{\widetilde{E}}\left(  K\right)  =\left(  \eta_{E}\left(  K\right)
+\eta_{L_{1}}\left(  K\right)  \right)  /2=\eta_{E}\left(  K\right)  /2
\]
because of $\eta_{L_{1}}\left(  K\right)  =0.$
\end{proof}

\section{References}

\begin{enumerate}
\item  Kadec M.I. and Pe\l czy\'{n}ski A. \textit{Bases, lacunary sequences
and complemented subspaces in the spaces }$L_{p}$ Studia Math. \textbf{21}
(1962) 161-176

\item  Aldous D.\textit{\ Subspaces of }$L_{1}$\textit{\ via random measures,}
Trans. AMS \textbf{267:2} (1981) 445-463

\item  Tokarev E.V. \textit{Reflexive subspaces of }$L_{1}$ (transl. from
Russian) Russian Mathematical Surveys.\textbf{40}, no.1 (1985) 247-248

\item  Semyonov E.M. \textit{Inclusion theorems for Banach spaces of
measurable functions }Dokl. A.N.SSSR \textbf{156} (1964) 1294-1295 (in Russian)

\item  Meyer P.A. Probability and potentials, \textit{N.Y. et al.: Blaisdell
Publ.Co,} 1966

\item  Aronszajn N. and Gagliardo E. \textit{Interpolation spaces and
interpolation methods}, Ann. Math. Pura App. Ser. \textbf{4, 68} (1965) 51-117

\item  Schwartz L. \textit{Geometry and probability in Banach spaces}, Bull.
AMS \textbf{4:2} (1981) 135-141

\item  Rodin V.A. and Semyonov E.M. \textit{Rademacher series in symmetric
spaces}, Analysis Math. \textbf{1:4} (1975) 207-222

\item  Tokarev E.V.\textit{\ Subspaces of symmetric spaces of functions
}(transl. from Russian)\ Funct. Anal. Appl. \textbf{13 }(1979) 152-153

\item  Tokarev E.V. \textit{Quotient spaces of Banach lattices and
Marcinkiewicz spaces} (transl. from Russian) Siberian Mathematical
Journal.\textbf{25} (1984) 332-338

\item  Novikov S.Ja., Semyonov E.M. and Tokarev E.V. \textit{On the structure
of subspaces of the spaces }$\Lambda_{p}\left(  \mu\right)  $ (transl. from
Russian) Translations. II. Ser. Amer. Math. Soc. \textbf{136} (1987) 121-127

\item  Sedaev A.A. \textit{On one Tokarev's theorem} Deposed in VINITI
11.07.78 no.\textbf{\ 2951-78}, 1-11 (in Russian)

\item  Lindenstrauss J. and Tzafriri L. \textit{Classical Banach spaces,}
Lecture Notes in Math. \textbf{338} (1973) 1-242

\item  Rudin W. \textit{Trigonometric series with gaps}, J. Math. Mech.
\textbf{9} (1960) 203-227

\item  Bennet G., Dor L.E., Gudman V., Johnson W.B. and Newmann C.M.
\textit{On uncomplemented subspaces of }\ $L_{p}$, $1<p<2$, Israel J. Math.
\textbf{26} (1977) 178-187
\end{enumerate}
\end{document}